\documentclass[11pt]{amsart}
\usepackage[T1]{fontenc}
\usepackage[ansinew]{inputenc}
\usepackage{amssymb}
\newtheorem{thm}{Theorem}[section]

\newtheorem{corollary}[thm]{Corollary}
\newtheorem{lem}[thm]{Lemma}

\newtheorem{example}{Example}
\newtheorem{proposition}[thm]{Proposition}
\newtheorem{definition}[thm]{Definition}
\newtheorem{remark}{Remark}[section]
\newcommand{\eps}{\varepsilon}
\newcommand{\const}{\mbox{const}}

\newcommand{\ad}{\mbox{ad}}
\newcommand{\divv}{\mbox{div}}
\newcommand{\curll}{\mbox{curl} }
\newcommand{\grad}{\mbox{grad}}
\newcommand{\spn}{\mbox{span}}

\newcommand{\DL}{\mbox{Det}_\mathcal{L}}

\begin{document}
\title{Solid Controllability in Fluid Dynamics}%
\author{Andrey A. Agrachev$^1$ \and Andrey
V. Sarychev$^2$}%
\address{$^1$International School for Advanced Studies (SISSA),
Trieste, Italy \& V.A.Steklov Mathematical Institute, Moscow, Russia
\newline  $^2$ DiMaD,
University of Florence, Italy}%
\email{agrachev@sissa.it,asarychev@unifi.it}%

\date{}%
%\dedicatory{}%
%\commby{}%
\begin{abstract}
We  survey results of recent activity towards studying
controllability and accessibility issues for equations of  dynamics
of incompressible fluids controlled by low-dimensional or,
degenerate, forcing. New results concerning controllability of
Navier-Stokes/Euler systems on two-dimensional sphere and on a
generic two-dimensional domain are represented.
\end{abstract}
% ----------------------------------------------------------------
\maketitle
% ----------------------------------------------------------------
%%%%%%%%%%%%%%%%%%%%%%%%%%%%%%%%%%%%%%%%%%%%%%%%%%%%%%%%%%%%%%%%%%%%

\vspace*{2mm}


\section{Introduction}
\label{intro} \markboth{A.A.Agrachev, A.V.Sarychev}{Solid
Controllability in Fluid Dynamics}

We survey results of recent activity  aimed at studying
controllability and accessibility properties of Navier-Stokes
systems (NS systems) controlled by low-dimensional (degenerate)
forcing. This choice of control is characteristic feature of our
problem setting. The corresponding equations  are
\begin{eqnarray}\label{eul1}
\partial u/\partial t + \nabla_u u + \grad p = \nu \Delta u +
F(t,x),  \\
\divv u=0.\label{ns2}
\end{eqnarray}
The words "{\it degenerate} forcing" mean that $F(t,x)$ can be
represented as:
$$F(t,x)=\sum_{k \in {\mathcal
K}^1}v_k(t)F^k(x), \ {\mathcal K}^1 \ \mbox{is finite}.$$ The word
"controlled" means that the functions $v_k(t), \ t \in [0,T]$
entering the forcing can be chosen freely among measurable
essentially bounded functions. In fact any functional space, which
is dense in $L_1[0,T]$ would fit.

The domains treated here include  2-dimensional (compact) Riemannian
manifolds $M$ homeomorphic to either sphere, torus or disc. The
latter case includes a rectangle and 2D simply connected domain $M$
with analytic boundary $\partial M$. We impose so called Lions
boundary condition, whenever boundary is nonempty.

Our approach stems from geometric control theory, which is
essentially based on differential geometry and Lie theory; geometric
control approach proved its effectiveness in studying controlled
dynamics in finite dimensions. We would like to report on some ideas
of how such methods can be extended onto the area of
infinite-dimensional dynamics and of controlled PDE. Extensions of
the geometric control theory onto infinite-dimensional case are
almost unknown. Classical Lie techniques are not well adapted for
infinite-dimensional case, and several analytic problems are
encountered.

In this contribution we concentrate almost exclusively on geometric
and Lie algebraic ideas of the accomplished work. For details on
analytic part we refer  interested readers to the publications
\cite{AS06,AS43,ARS,Rod0,Rod1,Shi1,Shi2}.

Applications of geometric theory to  the study of   controllability
of finite-di\-men\-sional systems is well established subject,
although many problems still remain unsolved. Starting point of the
activity aimed at controlling NS systems by degenerate forcing was
study (\cite{EM,ASpb,AS43,Rom}) of accessibility and controllability
of their finite-dimensional Galerkin approximations on
$\mathbb{T}^2$ and  $\mathbb{T}^2$ (periodic boundary conditions).
This question being settled, one should note, that controllability
of finite-dimensional Galerkin approximations of NS systems on many
other domains remains an open problem; answers for generic analytic
2D domains follow from results of Section~\ref{generdom}.

Study of infinite-dimensional case started in
\cite{ASDAN,AS43,AS06}, where we dealt with 2D NS/Euler system on 2D
torus $\mathbb{T}^2$ . In those publications notions of solid
controllability in projections and of approximate controllability
have been introduced and sufficient criteria for them have been
established.

To arrive to such criteria the technique of so-called Lie extensions
in infinite dimensions has been suggested. In the context of our
problem  this technique can be loosely interpreted as designing the
propagation to
 higher modes of the energy pumped by controlled forcing into lower
modes.

The control functions involved are fast-oscillating and analytic
part of the study consists of establishing continuity properties of
solutions of NS systems with respect to  so called relaxation metric
of forcing. The latter metric is weaker than the classical metrics
and is adapted for dealing with fast oscillating functions.

An extension of the above mentioned techniques techniques onto the
case of  NS system, subject to Lions boundary conditions on a
rectangle, has been accomplished by S.Rodrigues (\cite{Rod0}). In
the course of this study both geometric and analytic part needed to
be adjusted: Lie extensions turn more intricate and continuity
properties need to be reproved. These results are surveyed in
Section~\ref{rect}.

An additional analytic effort (perturbation techniques) is needed
for establishing controllability criteria for NS system on
smoothened rectangular; this is described in Section~\ref{smooth}.

A new approach is suggested for establishing  controllability  on a
generic 2D domain (Section~\ref{generdom}).

Finally study of Lie algebraic properties of spherical harmonics
results in controllability criterion for  NS/Euler system on 2D
sphere (Section~\ref{sphere}).

The results appearing in Sections~\ref{smooth}-\ref{sphere} have not
been previously published.

An interesting extension of the above described methods onto the
case of NS system on 3D torus has been accomplished by A.Shirikyan
in \cite{Shi1,Shi2}. The geometric part of his study coincides
essentially with the one of \cite{AS43} and of \cite{Rom}, but
additional  analytic difficulties in the 3D case are numerous. We do
not survey these results for the sake of remaining more geometric in
spirit and avoiding  new notation. Interested readers should consult
\cite{Shi1,Shi2}.

There was an extensive study of controllability of the Navier-Stokes
and Euler equations in particular by means of boundary control.
There are various results on exact local controllability of 2D and
3D Navier-Stokes equations obtained by A.Fursikov, O.Imanuilov,
global exact controllability for 2D Euler equation obtained by J.-M.
Coron, global exact controllability for 2D Navier-Stokes equation by
A.V.~Fursikov and J.-M.~Coron. The readers may turn to the book
\cite{Fu} and to the surveys \cite{FI} and \cite{Co} for further
references.

The authors are grateful to S.Rodrigues for useful comments and help
during preparation of this contribution.

\section{2D NS/Euler system controlled by degenerate forcing.
Definitions and problem setting} \label{dset}

\subsection{NS/Euler system on 2D Riemannian manifold}

Representation of  NS/Euler system in the form
(\ref{eul1})-(\ref{ns2}) requires interpretation, whenever one
considers the system on 2D domain $M$ with arbitrary Riemannian
metric. There is a general  way of representing NS/Euler systems on
any $n$-dimensional Riemannian manifold (see e.g. \cite{ArnKh}), but
we prefer  to remain in $2$ dimensions and to advance with some
elementary vector  analysis in the 2D Riemannian case.

We consider smooth (or analytic) 2-dimensional Riemannian manifold
$M$ (with or without boundary), endowed with Riemannian metric
$(\cdot , \cdot)$ and with area 2-form $\sigma$. All functions,
vector fields, forms, we deal with, will be smooth.

Any vector field $y$ on $M$ can be paired  with two differential
1-forms defined as:
\begin{eqnarray}
  y \mapsto y^\flat: \ \langle y^\flat, \xi \rangle= (y, \xi) \nonumber \\
y \mapsto y^\sharp: \ \langle y^\sharp, \xi \rangle  = \sigma (y ,
\xi),\label{flana}
\end{eqnarray}
for each vector field $\xi$.  Obviously $\langle y^\sharp , y
\rangle = \sigma (y , y)=0$.

Note that for any $1$-form $\lambda$ there holds:
\begin{equation}\label{nat}
 \lambda \wedge y^\sharp = \langle \lambda , y \rangle \sigma .
\end{equation}
To prove (\ref{nat}) it suffices to compare the values of $2$-forms
$\lambda \wedge y^\sharp$ and $ \langle \lambda , y \rangle \sigma $
on any pair of linearly independent vectors. Evidently (\ref{nat})
is valid if $y$ (and $y^\sharp$) vanishes. If $y \neq 0$, we take a
pair $y,z$, which is linearly independent. Then
$$(\lambda \wedge y^\sharp)(y,z)=\left|
                                     \begin{array}{cc}
                                     \langle \lambda , y \rangle & \langle y^\sharp , y\rangle \\
                                      \langle \lambda ,z \rangle & \langle y^\sharp , z \rangle \\
                                     \end{array}
                                   \right|=\langle \lambda , y \rangle \sigma(y,z).  $$

Now we define the vorticity $\curll$ and the divergence $\divv$ of a
vector field via  the differentials $d y^\flat , d y^\sharp$. These
latter are 2-forms; we put
$$d y^\flat =(\curll y)\sigma , d y^\sharp=(\divv y)\sigma , $$
or by abuse of notation:
\begin{equation}\label{cudi}
(\curll y)=d y^\flat /\sigma , (\divv y)=d y^\sharp/\sigma .
\end{equation}

The gradient $\grad \varphi$ of a function $\varphi$ is the vector
field paired with $d \varphi$ metrically: $(\grad \varphi)^\flat=d
\varphi$.

As in the Euclidean case the vorticity  of gradient vector field
vanishes:
$$\curll (\grad \varphi)=d (\grad \varphi)^\flat/\sigma =d (d
\varphi)/\sigma=0.$$

While in 3D case operator $\curll$ transforms vector fields into
vector fields,  in the  2D case it transforms vector fields into
scalar functions (actually in a component of vector field directed
along additional third dimension). We will define the vorticity
operator $\curll$ on functions. The result of its action on a
function $\phi$ is  a vector field $\curll \phi$, which satisfies
the relation:
$$\langle \lambda ,\curll \phi \rangle \sigma=(d\phi \wedge \lambda),$$
for each $1$-form $\lambda$. By virtue of (\ref{nat}) and due to the
nondegeneracy of paring $y \mapsto y^\sharp$ we conclude:
\begin{equation}\label{curnat}
(\curll \phi)^\sharp=- d \phi.
\end{equation}

As in the Euclidean case the divergence of  vorticity of a function
vanishes:
\begin{equation}\label{dicu}
\divv (\curll \phi)= d(\curll \phi)^\sharp /\sigma=-d (d \phi)
/\sigma =0.
\end{equation}

Coming back to the equation (\ref{ns2})  we note that the condition
$\divv u=0$ can be written down as
\begin{equation}\label{du0}
 d u^\sharp=0.
\end{equation}
 If
$M$ is simply connected we conclude that $u^\sharp$ must be a
differential: $u^\sharp=-d \psi$, where $\psi$ is so called {\it
stream function}. By virtue of (\ref{curnat}) $$\curll \psi=u.$$

For non simply connected domains we {\it impose}  a condition which
guarantees the exactness; in the next subsection we comment on it.

Given the symplectic structure on $M$ defined by $\sigma$ and
$(\cdot, \cdot)$ we see that $u$ is Hamiltonian vector field
corresponding to the Hamiltonian $-\psi$:
$u=-\overrightarrow{\psi}$.

The nonlinear term $\nabla_u u$ in the right-hand side of
(\ref{eul1}) corresponds to covariant derivative of the Riemannian
(metric torsion-free) connection on $M$.

Finally we define the Laplace-Beltrami operator $\Delta$ as
$$\Delta=\curll^2.$$ In Hodge theory (see \cite{ArnKh}) this operator
transforms $p$-forms into $p$-forms; in our notation $\Delta$
transforms  vector fields into vector fields and functions into
functions.

\subsection{Helmholtz form of 2D NS sytem}
To arrive to the Helmholtz form of the NS system
(\ref{eul1})-(\ref{ns2}), we  apply the operator $\curll$ to both
parts of (\ref{eul1}). As a result we get for the vorticity $\curll
u=w$ the equation
\begin{equation}\label{eul13}
   \partial
w/\partial t  + \curll (\nabla_u u)  = \nu \Delta w + f(t,x),
\end{equation}
where $f(t,x)=\curll F(t,x)$.

One should note that the vorticity of  $\grad p$ vanishes and that
the operator $\curll$ commutes with $\Delta = \curll^2$ .

To calculate $\curll (\nabla_u u)$ according to the formula
(\ref{cudi})we first compute  the $1$-form $(\nabla_u u)^\flat$,
adapting the argument of
 \cite[\S IV.1.D]{ArnKh}.

 Let $y$ be a vector field which commutes with $u$: the Lie-Poisson
bracket $[u,y]=0$. Then
\begin{equation}\label{flat}
\langle (\nabla_u u)^\flat , y \rangle =(\nabla_u u , y)=L_u(u, y)-
(u,\nabla_u y).
\end{equation}
 (Here and below  $L_u$ denotes Lie derivative. Note that for the covariant derivative of metric
connection there holds: $L_u(u, y)=(\nabla_u u, y)+(u,\nabla_u y)$.)
Since the connection is torsion-free and $[u,y]=0$, then $\nabla_u
y-\nabla_y u=0$, and the right-hand side of  (\ref{flat}) can be
represented as
$$L_u \langle u^\flat , y \rangle - (u,\nabla_y u)=L_u \langle u^\flat , y \rangle -
\frac{1}{2}\langle d(u,u), y \rangle .$$ Besides $L_u \langle
u^\flat , y \rangle= \langle L_u u^\flat , y \rangle,$ as long as
$L_u y=[u,y]=0$, and we obtain
$$\langle (\nabla_u u)^\flat , y \rangle =\langle L_u u^\flat , y \rangle -
\frac{1}{2}\langle d(u,u), y \rangle .$$ As far as one can find
vector field $y$, which commutes with $u$ and has any prescribed
value at a given point, we conclude: $(\nabla_u u)^\flat =L_u
u^\flat - \frac{1}{2} d(u,u)$.

Using the definition of $\curll$ (\ref{cudi}) we get
$$\curll (\nabla_u u)= d((\nabla_u u)^\flat)/\sigma= d L_u  u^\flat/\sigma=
L_u d u^\flat/\sigma=L_u (w \sigma)/ \sigma= L_u w .$$ Hence $
\curll (\nabla_u u)=L_u w$.

For $u$ being  Hamiltonian vector field with the Hamiltonian $-\psi$
there holds $\nabla_u w=-\{\psi ,w\} $, where $\{\cdot , \cdot\}$ is
Poisson bracket of functions.

The Helmholtz form of the Navier-Stokes equation (cf. \cite{ArnKh})
reads
  $$   \frac{\partial w}{\partial t}-\{\psi,w\}-\nu\Delta
w=f(t,x).$$

Note that $w=\curll u=\curll^2 \psi=\Delta \psi$.

The Lions  condition written in terms of the vorticity $w$ and the
stream function $\psi$ reads:
\begin{equation}\label{bc}
    \psi\bigr|_{\partial M}=w\bigr|_{\partial M}=0.
\end{equation}
If the boundary $\partial M$ of $M$ is smooth, then the Hamiltonian
vector field $u=-\overrightarrow{\psi}$ is tangent to $\partial M$.

Given the vorticity $w$ and  the boundary conditions (\ref{bc}), one
can recover in a unique way the velocity field $u$, which
corresponds to an {\it exact}  1-form $u^\sharp$. The corresponding
formula is $u=\curll \psi$, where $\psi$ is the unique solution of
the Dirichlet problem $\Delta \psi =w$ under boundary condition
(\ref{bc}). Indeed such $u$ is divergence free  and its vorticity
equals $w$ by the definition of $\Delta$.

The NS system can be written as
\begin{equation}\label{evort}
    \frac{\partial w}{\partial t}-\{\Delta^{-1}w,w\}-\nu\Delta
w=f(t,x).
\end{equation}

This last equation looks universal: in fact its dependence on the
domain is encoded in the properties of the Laplacian $\Delta$ on
this domain. It is well explained in \cite{Arn1,ArnKh} that the
Euler equation for fluid motion is infinite-dimensional analogy of
Euler equation for rotation of (multidimensional) rigid body, and
the Laplacian in (\ref{evort}) plays role of tensor of inertia of
rotating rigid body.

\subsubsection{Stream function on flat torus}
Let us consider  the flat torus $\mathbb{T}^2$, endowed with
standard Riemannian metric and with area form $\sigma$,  both
inherited from covering of $\mathbb{T}^2$ by the Euclidean plane.
Let  $\varphi_1, \varphi_2$ be the 'euclidean' coordinates on
$\mathbb{T}^2$. We proceed in the space of velocities $u$ with
vanishing space average: $\int_{\mathbb{T}^2}u d\sigma =0$ (due to
flatness we may think that all the velocities belong to the same
linear space).

   To establish exactness of the closed
$1$-form $u^\sharp$  (involved in  (\ref{du0})) it suffices to prove
that its integral along a generator of a torus vanishes.  By Stokes
theorem the integrals of the closed form $u^\sharp$ along any two
homologous paths have the same value.

 Taking $u=(u_1,u_2)$ we get $u^\sharp=-u_2d
\varphi_1+u_1d \varphi_2$. Integrating $u^\sharp$ along a loop
$\Gamma: \ \varphi_1=\alpha$, we obtain the value of the integral
$\int_{\Gamma}u^\sharp=\int_0^{2\pi}u_1 d \varphi_2 =c(\alpha)$,
which by the aforesaid is constant: $c(\alpha) \equiv c$.
Integrating it with respect to $\varphi_1$ we conclude $2\pi
c=\int_0^{2\pi}u_1 d \varphi_2 d \varphi_1=\int_{\mathbb{T}^2}u_1
d\sigma =0$. Hence $\int_{\Gamma}u^\sharp=c=0$. The same holds for
the loops $\Gamma': \ \varphi_2=\const$.

\subsection{Controllability: definitions}
In what follows we reason in terms of so called {\it modes} which
are the eigenfunctions $\phi^k(x)$ of the Laplace-Beltrami operator
$\Delta$ defined in the space of vorticities $w$: $\Delta
\phi^k(x)=\lambda_k \phi^k(x)$.

Representing  $w,f$ in (\ref{evort}) as a series $w(t,x)=\sum_k
q_k(t)\phi^k(x),\ f(t,x)=\sum_k v_k(t)\phi^k(x)$ with respect to the
basis of eigenfunctions one can write  the NS system as an infinite
system of ODE on the coefficients $q_k(t)$. Assume
$$\{\phi^i(x), \phi^j(x)\} =\sum_{k}C^{ij}_k\phi^k(x).$$
Then the equation (\ref{evort}) can be written down 'in  coordinate
form' as
\begin{equation}\label{coord}
    \dot{q}_k-\sum_{i,j}C^{ij}_k \lambda^{-1}_iq_iq_j-\nu \lambda_k q_k=v_k(t).
\end{equation}

Typically we will consider the controlled forcing, which is applied
to few modes
 $\phi^k(x), \ k \in \mathcal{K}^1, \ \mathcal{K}^1$ - finite. Then  in
the system (\ref{coord})
 the controls
enter only the equations indexed by $k \in \mathcal{K}^1$, while for
$k \not\in \mathcal{K}^1, \ v_k=0$.

 Introduce another finite set $\mathcal{K}^{o}$ of {\em observed} modes.
 We will always assume $\mathcal{K}^{o} \supset
\mathcal{K}^1$.  We identify the space of observed modes with
$\mathbb{R}^N$ and denote by $\Pi^{o}$ the operator of projection of
solutions onto the space of observed modes $\spn\{\phi_k| \ k \in
\mathcal{K}^{o}\}$. The coordinates corresponding to the observed
modes are reunited in {\it observed component} $q^{o}$.

 Galerkin
$\mathcal{K}^{o}$-approximation of the 2D NS/Euler system is the
equation (ODE) for $q^{o}(t)$, obtained by projecting the 2D NS
system onto the space of observed modes and putting all the
components $q_k(t), \  k \not\in \mathcal{K}^{o}$ zero. The
resulting equation is
\begin{equation}\label{galobs}
     \frac{\partial q^{o}}{\partial t}-\Pi^{o}\{\Delta^{-1}q^{o},q^{o}\}-\nu\Delta
q^{o}=f(t,x).
\end{equation}
As far as  $\mathcal{K}^{o} \supset \mathcal{K}^1$, i.e. controlled
forcing $f$ only affects part of observed modes, then
$\Pi^{o}f(t,x)=f(t,x)$.

In coordinate form passing to Galerkin approximation means omitting
 the equations  (\ref{coord}) for variables  $q_k$ with  $k
\not\in \mathcal{K}^{o}$ and taking these $q_k$ equal $0$ in the
resting equations.

We say that a control $f(t,x)$ {\it steers the system (\ref{evort})
(or (\ref{galobs})) from $\tilde{\varphi}$ to $\hat{\varphi}$ in
time $T$}, if for equation (\ref{evort}) forced by $f$ the solution
with the initial condition $\tilde{\varphi}$ at $t=0$ takes 'value'
$\hat{\varphi}$ at $t=T$.

The first notion of controllability under study is controllability
of Galerkin approximation.

\begin{definition}[controllability of Galerkin approximation]
Galerkin $\mathcal{K}^{o}$-approximation of 2D  NS/Euler systems is
time-$T$ globally controllable if for any two points $\tilde{q},
\hat{q}$  in $\mathbb{R}^N$,  there exists a control which steers in
time $T$ this Galerkin approximation from $\tilde{q}$ to $\hat{q}. \
\Box$
\end{definition}
This is purely finite-dimensional notion. The next notion  regards
finite-dimensional component of solutions, but takes into account
complete infinite-dimensional dynamics. Let us introduce some
terminology.

%%\begin{definition}{\it (controllability in finite-dimensional
%%projection)} \label{coproj}
%% Let $\mathcal L$ be a finite-dimensional subspace
%%of $H_2(\mathbb{T}^2)$ and $\Pi^{\mathcal L}$ be
%% $L_2$-orthogonal projection of $H_2(\mathbb{T}^2)$ onto $\mathcal L$.
%%  The 2D NS/Euler  system is  time-$T$ globally
%%controllable in projection on $\mathcal L$ if for any
%%$\tilde{\varphi} \in H_2(\mathbb{T}^2)$ and for any
%% $\hat{q} \in \mathbb{R}^N$
%%there exists a control which steers the  system in time $T$ from
%%$\tilde{\varphi}$ to some $\hat{\varphi} \in (\Pi^{\mathcal
%%L})^{-1}(\hat{q}). \ \Box$
%%\end{definition}

%%\begin{definition}{\it ($L_2$-approximate controllability)}
%%\label{dapc}
%% The 2D NS/Euler system is  time-$T$
%%$L_2$-appro\-xi\-mate\-ly controllable, if for any two points
%%$\tilde{\varphi}, \hat{\varphi} \in H_2$ and for any $\eps >0$ there
%%exists a control which steers the  system in time $T$ from
%%$\tilde{\varphi}$ to the $\eps$-neighborhood of $\hat{\varphi}$ in
%%$L_2$-norm. $\Box$
%%\end{definition}

\begin{definition}[attainable sets of NS systems]
Attainable set $\mathcal{A}_{\tilde{\varphi}}$   of the  NS/Euler
system (\ref{evort}) is the set of points in $H^2(M)$ attained from
$\tilde{\varphi}$ by means of essentially bounded measurable
controls in any positive time. For each $T>0$ time-$T$ (time-$\leq
T$)  attainable set $\mathcal{A}^T_{\tilde{\varphi}}$
($\mathcal{A}^{\leq T}_{\tilde{\varphi}}$) of the  NS/Euler system
is the set of points attained from $\tilde{\varphi}$ by means of
essentially bounded measurable controls in time $T$ (in time $\leq
T$). Attainable set $\mathcal{A}_{\tilde{\varphi}}=\bigcup_T
\mathcal{A}^{T}_{\tilde{\varphi}}. \ \Box$
\end{definition}

\begin{definition}\label{dapc}
The NS/Euler system is time-$T$ globally controllable in projection
onto $\mathcal L$ if  for each $\tilde{\varphi}$ the image
$\Pi^{\mathcal L}\left(\mathcal{A}^T_{\tilde{\varphi}}\right)$
coincides with ${\mathcal L}$.
\end{definition}

\begin{definition}\label{acon}
The NS/Euler system is time-$T$ $L_2$-approximately controllable if
$\mathcal{A}^T_{\tilde{\varphi}}$ is $L_2$-dense in $H^2. \ \Box$
\end{definition}

Let us introduce the notion of accessibility in projection.

\begin{definition}{\it (accessibility in finite-dimensional
projection)} \label{acproj}
 Let $\mathcal L$ be a finite-dimensional subspace of $H_2(M)$ and $\Pi^{\mathcal L}$ be
 $L_2$-orthogonal projection of $H_2(M)$ onto $\mathcal L$.
  The NS/Euler  system is  time-$T$ accessible in projection on $\mathcal L$ if for any
$\tilde{\varphi} \in H_2(M)$ the image  $\Pi^{\mathcal
L}\left(\mathcal{A}^T_{\tilde{\varphi}}\right)$ contains interior
points in ${\mathcal L}. \ \Box$
\end{definition}

Now we will introduce notion  of  {\it solid controllability}.

\begin{definition}\label{itomap}
Fix initial condition $\tilde{\varphi} \in H_2(M)$
 for trajectories of the controlled 2D NS/Euler system. Let $v(\cdot) \in
L_\infty\left([0,T];\mathbb{R}^r\right)$ be the controlled forcing
and $w_t$ be the corresponding trajectory of the NS system.

%%The correspondence between the controlled forcing $v(\cdot) \in
%%L_\infty\left([0,T];\mathbb{R}^d\right)$ and the corresponding
%%trajectory (solution) $w_t$ of the  system is established by
%%forcing/trajectory map ($\mathcal{F/T}$-map).

%%The correspondence between the controlled forcing $v(\cdot)$ and the
%%observed component $q(t)=\Pi^{o}w_t$ (an $\mathbb{R}^N$-valued
%%function) of the corresponding trajectory is established by
%%forcing/observation map ($\mathcal{F/O}$-map).

If NS/Euler system is considered on an interval $[0,T]\ (T <
+\infty)$, then the map $E_T : v(\cdot) \mapsto w_T$ is called
end-point map; the map $\Pi^{o} \circ \mathcal{F/T}_T$ is called
end-point component map, the composition $\Pi^{\mathcal L} \circ
\mathcal{F/T}_T$ is called $\mathcal L$-projected end-point map.
$\Box$
\end{definition}

%%\begin{remark}
%%In the terminology of control theory the first  two maps would be
%%called input/trajectory and input/output maps  correspondingly.
%%$\Box$
%%\end{remark}

%%\begin{remark}
%%  Evidently time-$T$ controllability of the NS/Euler system in
%%observed component or in finite-dimensional projection  is the same
%%as surjectiveness of the corresponding end-point maps. $\Box$
%%\end{remark}

\begin{definition}\label{socov}
Let $\Phi: \mathcal{M}^1 \mapsto \mathcal{M}^2$ be a continuous map
between two metric spaces, and $S \subseteq \mathcal{M}^2$ be any
subset. We say that $\Phi$ covers $S$ solidly, if $S \subseteq
\Phi(\mathcal{M}^1)$ and this inclusion is stable with respect to
$C^0$-small perturbations of $\Phi$, i.e. for some
$C^0$-neighborhood $\Omega$ of $\Phi$ and for each map $\Psi \in
\Omega$, there holds: $S \subseteq \Psi(\mathcal{M}^1). \ \Box$
\end{definition}

%%In what follows $M^2$ will be finite-dimensional vector space.

\begin{definition}{\it (solid controllability in finite-dimensional
projection)} \label{socoproj}
   The 2D NS/Euler  system is  time-$T$ solidly globally
controllable in projection on finite-dimensional subspace $\mathcal
L \subset H^2(M)$, if for any bounded set $S$  in $\mathcal L$ there
exists a set of controls $B_S$ such that $\left(\Pi^{\mathcal L}
\circ \mathcal{F/T}_T\right)(B_S)$ covers $S$ solidly. $\Box$
\end{definition}

\subsection{Problem setting}

In this contribution we will discuss  the following questions.

\begin{itemize}
  \item Under what conditions the 2D NS/Euler system  is
globally controllable in observed component? $\Box$
  \item Under what conditions the 2D NS/Euler system  is
solidly controllable in a finite-dimensional projection? $\Box$
  \item Under what conditions the 2D NS/Euler system  is
accessible in a finite-dimensional projection? $\Box$
  \item Under what conditions the 2D NS/Euler system is
$L_2$-appro\-xima\-tely controllable? $\Box$
\end{itemize}

As we explained above the geometry of controllability is encoded in
spectral properties of the Laplacian $\Delta$ and therefore on the
geometry of  the domain on which the controlled NS system evolves.
Below we provide answers for particular types of domains.

\section{Geometric control: accessibility and controllability via Lie brackets}
\label{tools}

In this section we collect some results of geometric control theory
regarding accessibility and controllability of  {\it
finite-dimensional} real-analytic control-affine systems of the form
\begin{equation}\label{caf}
\dot{x}=f^0(x)+ \sum_{i=1}^r f^i(x)v_i(t), \ x(0)=x^0, \ v_i(t) \in
\mathbb{R}, \ i=1, \ldots ,r.
\end{equation}
Geometric approach is coordinate-free, so that it is adapted for
dealing with dynamics on manifolds, but we will assume that  the
system (\ref{caf}) is defined on a finite-dimensional linear space
$\mathbb{R}^N$ in order to maintain parallelism
 with NS systems, which evolve in Hilbert spaces.

We use standard notation $P_t=e^{tf}$ for the flow corresponding to
a vector field $f$.

\subsection{Orbits, Lie rank  accessibility}

%%%\begin{definition}[attainable sets]
Let $v(\cdot) \in L_\infty\left([0,T];\mathbb{R}^r\right)$ be
admissible controls and $x(t)$ be corresponding trajectories of the
 system $\dot{x}=f^0(x)+ \sum_{i=1}^r f^i(x)v_i(t)$ with initial point $x(0)=x^0$.
We again introduce the end-point map $\mathcal{E}_T: v(\cdot)
\mapsto x_v(T)$; here $x_v(\cdot)$ is the trajectory of (\ref{caf})
corresponding to the control $v(\cdot)$.

For each $T>0$ time-$T$ (time-$\leq T$)  attainable set
$\mathcal{A}^T_{x^0}$  from $x^0$
 of the system (\ref{caf}) is the image of the set
$L_\infty\left([0,T];\mathbb{R}^r\right)$ under the map
$\mathcal{E}_T$, or, equivalently,  the set of  points $x(T)$
attained in time $T$ from $x^0$ by means of admissible controls. The
time-$\leq T$  attainable set from $x^0$ is $\mathcal{A}^{\leq
T}_{x^0} \bigcup_{t \in [0,T]} \mathcal{A}^{t}_{x^0}$.  Attainable
set from $x^0$ of the system (\ref{caf}) is
$\mathcal{A}_{x^0}=\bigcup_{T \geq 0} \mathcal{A}^{T}_{x^0}. \ \Box$
%%%\end{definition}

Important notions  of geometric control theory are orbits of control
system.

\begin{definition}[orbits and zero-time orbits of control systems] An
orbit of the control system (\ref{caf})  passing through  $x^0$  is
the set of points obtained from $x^0$ under the action of (the group
of) diffeomorphisms  of the form $e^{t_1 f^{u^1}} \circ \cdots \circ
e^{t_N f^{u^N}}$, where $t_j \in \mathbb{R}, \ j=1, \ldots , N,$ and
$f^{u^j}=f^0+\sum_{i=1}^r f^i(x)u^j_i$ is the right-hand side of
(\ref{caf}) corresponding  to constant control $u^j =(u^j_1, \ldots
, u^j_r) \in \mathbb{R}^r$. Zero-time orbit is the subset of the
orbit, resulting from the action of these diffeomorphisms  subject
to condition $\sum_j t_j=0.\ \Box$
\end{definition}

If we consider 'symmetrization' of the system (\ref{caf}),
$$\dot{x}=f^0(x)v_0+ \sum_{i=1}^r f^i(x)v_i(t), \ x(0)=x^0, \ v_0 \in \mathbb{R}, \
v_i(t) \in \mathbb{R}, \ i=1, \ldots ,r,
 $$
 then the orbit of (\ref{caf}) can be interpreted as the attainable
set from $x^0$ of this symmetrization corresponding to application
of piecewise-constant controls.

The famous Nagano theorem relates properties of the orbits and Lie
algebraic properties of the system. It claims that the orbit and the
zero-time orbit of an analytic system (\ref{caf}) are immersed
manifolds of $\mathbb{R}^N$, and tangent spaces to these orbits can
be calculated via Lie brackets of the vector fields $\{f^0, \ldots
,f^m\}$.

\begin{definition}[Lie rank and zero-time Lie rank]
\label{lrank} Take  the Lie algebra $\mbox{Lie}\{f^0, \ldots ,
f^m\}$ generated by $\{f^0, \ldots , f^m\}$  and evaluate vector
fields from  $\mbox{Lie}\{f^0, \ldots , f^m\}$ at a point $x$; the
dimension of the resulting linear space $\mbox{Lie}_x\{f^0, \ldots ,
f^m\}$ is Lie rank  of the system $\{f^0, \ldots , f^m\}$ at $x$.

Take the Lie ideal generated by  $\spn\{f^1, \ldots , f^m\}$ in
$\mbox{Lie}\{f^0, \ldots , f^m\}$ and evaluate vector fields from it
at $x$;  the dimension of the resulting linear space
$\mbox{Lie}^0_x\{f^0, \ldots , f^m\}$,  is zero-time Lie rank at $x$
of the system $\{f^0, \ldots , f^m\}. \ \Box$
\end{definition}

These two Lie ranks either are equal or differ by $1$.

The Nagano theorem claims that in analytic case $\mbox{Lie}_x\{f^0,
\ldots , f^m\}$ and $\mbox{Lie}^0_x\{f^0, \ldots , f^m\}$ are
tangent spaces at each point $x$ of the orbit and zero-time orbit
respectively.

Accessibility properties of analytic control system (\ref{caf}) are
determined by its Lie ranks. Recall that the system is called
accessible if its attainable set $\mathcal{A}_{x^0}$ has nonempty
interior and is strongly accessible if $\forall T>0$ attainable sets
$\mathcal{A}^T_{x^0}$ have nonempty interior. The following fact
holds.

\begin{thm}[Jurdjevic-Sussmann ($C^\omega$ case),Krener ($C^\infty$ case)]
\label{kren} If the Lie rank of the system of vector fields $\{f^0,
\ldots , f^r\} $ at $x^0$ equals $n$, then $\forall T>0$ the
interior of the attainable set $\mathcal{A}^{\leq T}_{x^0}$ is
nonvoid. Besides $\mathcal{A}^{\leq T}_{x^0}$ possesses the interior
which is dense in it. If zero-time Lie rank at $x^0$ equals $n$,
then $\forall T>0$ the interior of the attainable set
$\mathcal{A}^{T}_{x^0}$ is nonvoid and this interior is dense in
$\mathcal{A}^{T}_{x^0}. \ \Box$
\end{thm}

See \cite{Ju,ASkv} for the proof.

Let ${\mathcal L}$ be a linear subspace of $\mathbb{R}^N$ and
$\Pi^{\mathcal L}$ be a projection of $\mathbb{R}^N$ onto ${\mathcal
L}$. A control system (\ref{caf}) is accessible (strongly
accessible) from $x$ in projection on $\mathcal L$ if the image
$\Pi^{\mathcal L} \mathcal{A}_{x^0}$ ($\Pi^{\mathcal L}
\mathcal{A}^T_{x^0}$) contains interior points in ${\mathcal L}$
(for each $T>0$).

One easily derives from the previous theorem a criterion for
accessibility in  projection.

\begin{thm}
If $\Pi^{\mathcal L}$ maps $\mbox{Lie}_{x}\{f^0, \ldots , f^m\}$
(respectively $\mbox{Lie}^0_x\{f^0, \ldots , f^m\}$) onto ${\mathcal
L}$, then control system (\ref{caf}) is accessible (respectively
strongly accessible) at $x$ in  projection on $\mathcal L. \ \Box$
\end{thm}

\begin{proof}
The proofs of the two statements are similar; we  sketch the proof
of the first one. Consider the orbit of the system (\ref{caf})
passing through $x_0$. The tangent space to the orbit at each of its
points $x$ coincides with $\mbox{Lie}_{x}\{f^0, \ldots , f^m\}$.

By  Theorem~\ref{kren} the attainable set of the system possesses
relative interior with respect to the orbit. Moreover there are
interior points $x_{int} \in \mathcal{A}_{x^0}$ arbitrarily close to
$x^0$ so that $\Pi^{\mathcal L}$ maps $\mbox{Lie}_{x_{int}}\{f^0,
\ldots , f^m\}$ onto ${\mathcal L}$. Then sufficiently small
neighborhoods of $x_{int}$ in the orbit are contained in
$\mathcal{A}_{x^0}$ and are mapped by $\Pi^{\mathcal L}$ onto a
subset of ${\mathcal L}$ with nonempty interior.
\end{proof}

\subsection{Lie extensions and controllability}
Controllability is  stronger and much more delicate property than
accessibility.  To verify it it does not suffice in general to
compute the Lie rank which accounts for {\it all} the Lie brackets.
Instead one should select  'good Lie brackets' avoiding 'bad Lie
brackets' or 'obstructions'.

To have a general idea of what good and bad Lie brackets can be like
let us consider the following elementary example.

\begin{example}\label{ex1}
$$\dot{x}_1=v, \ \dot{x}_2=x_1^2.$$
\end{example}
This is $2$-dimensional control-affine system (\ref{caf}) with
$f^0=x_1^2 \partial / \partial x_2, \ f^1=\partial / \partial x_1$.
The Lie rank of this system equals $2$ at each point, the system is
accessible but uncontrollable from each point $\hat{x}=(\hat{x}_1,
\hat{x}_2)$ given the fact that we can not achieve points with $x_2
< \hat{x}_2$. One can prove that the attainable set
$\mathcal{A}_{\hat{x}}$ coincides with the half-plane $x_2 >
\hat{x}_2$ with added point $\hat{x}$. One sees that that it is
possible to move freely (bidirectionally) along good vector field
$f^1$, while along bad Lie bracket $[f^1[f^1,f^0]]= 2\partial /
\partial x_2$ we can move only in one direction.

%%%%%%%%%%%%%%%%%%%%%%%%%%%%%%%%%%%%%%%%%%%%%%%%%%%%%%%%%%%%%%%%

 Good Lie brackets form  {\em Lie extension} of a control system.

\begin{definition}\label{ext}
The family $\mathcal{F}'$   of real analytic vector fields is:

i) an extension of $\mathcal{F}$ if $\mathcal{F}' \supset
\mathcal{F}$ and the closures of the attainable sets
$\mathcal{A}_{\mathcal{F}}(\tilde{x})$ and
$\mathcal{A}_{\mathcal{F}'}(\tilde{x})$ coincide;

ii) time-T extension of $\mathcal{F}$ if $\mathcal{F}' \supset
\mathcal{F}$ and $\forall T
>0$ the closures of the time-$T$ attainable sets
$\mathcal{A}^T_{\mathcal{F}}(\tilde{x})$ and
$\mathcal{A}^T_{\mathcal{F}'}(\tilde{x})$ coincide;

iii) fixed-time extension if it is time-$T$ extension $\forall T>0$.

The vector fields from $\mathcal{F}'\setminus  \mathcal{F}$ are
called i) compatible, ii) compatible in time $T$,  iii) compatible
in a fixed time  with $\mathcal{F}$ in the cases i), ii) and iii)
correspondingly.$\ \Box$
\end{definition}

The inclusions $\mathcal{A}_{\mathcal{F}}(\tilde{x}) \subset
\mathcal{A}_{\mathcal{F}'}(\tilde{x}), \
\mathcal{A}^T_{\mathcal{F}}(\tilde{x}) \subset
\mathcal{A}^T_{\mathcal{F}'}(\tilde{x})$ are obvious.  Less obvious
is the following Proposition (see \cite{ASkv}).

\begin{proposition}
\label{incl} If an extension $\mathcal F'$ of an analytic system
$\mathcal F$ is globally controllable, then $\mathcal F'$ also is.
$\Box$
\end{proposition}

\begin{remark}
When talking about time-$T$ extensions one can consider also
extensions by time-variant vector fields $X_t, \ t \in [0,T]$. We
will say that vector field $X_t$ is time-$T$ compatible with
$\mathcal{F}$ if it drives the system in time $T$ from $\tilde{x}$
to the closure of $\mathcal{A}^T_{\mathcal{F}}(\tilde{x}). \ \Box$
\end{remark}

Our idea is to proceed with a series of extensions of a control
system in order to arrive to extended system for which the
controllability can be verified and then apply
Proposition~\ref{incl}.

Obviously the Definition~\ref{ext} is nonconstructive. In what
follows we will use three  particular types of extensions.

The first natural type is based on possibility of taking topological
closure of a set of vector fields, maintaining closures of
attainable sets.

\begin{proposition}(see \cite[Ch. 3, \S2, Th.5]{Ju})
Topological (with respect to $C^\infty$ convergence on compact sets)
closure $cl(\mathcal F)$ of $\mathcal F$ is Lie extension. $\Box$
\end{proposition}

The second method  underlies the theory of relaxed (or sliding mode)
controls started \cite{Gam,G60}, which provided an extension of
pioneering contributions by L.C.Young \cite{Yng} and E.J. McShane
\cite{McSh} onto the field of optimal control theory. To introduce
it we consider a family of so-called relaxation seminorms
$\|\cdot\|_{s,K}$ of time-variant vector fields $X_t, \ t \in
[0,T]$:
\begin{equation}\label{reno}
  \|X_\cdot\|^\mathrm{rx}_{s,K}=\max_{t\in[0,T]}\left|
\int\limits_0^t\|X_\tau\|_{s,K} d\tau\right|,
\end{equation}
with $K$ being a compact in $\mathbb{R}^N$, $s \geq 0$ being an
integer and $\|X_\tau\|_{s,K}$ being a $C^k$-norm on the compact
$K$. The family of relaxation seminorms define relaxation
($C^\infty$) topology (metrics) in the set of time-variant vector
fields.

\begin{proposition}[see \cite{Gam, ASkv}]\label{fo1}
Let a sequence  of time-variant vector fields $X^j_t$ converge to
$X_t$ in the relaxation metrics and let these vector fields have
compact supports. Then the flows of the vector fields $X^j_t$
converge to the flow of the vector field $X_t. \ \Box$
\end{proposition}

As a corollary of this result one can prove

\begin{proposition}
 For the systems $$\mbox{co}\mathcal{F}=\left\{\sum_{i=1}^m \beta_i f_i, \ f_i \in
\mathcal{F}, \ \beta_i \in C^\omega(\mathbb{R}^N), \ \beta_i \geq 0,
\ \sum_{i=1}^m \beta_i \equiv 1, \ i=1, \ldots , m\right\},$$ and
$\mathcal{F}$ the closures of their time-$T$ attainable sets
coincide. $\Box$
\end{proposition}

Proof of this result and of its modifications  can be found in
\cite[Chapter 8]{ASkv}, \cite[Chapter 3]{Ju},\cite[Chapters II,III
]{Gam}.

 The next type of
extension, we will  use, relies upon  Lie brackets. It appeared in
our earlier work on controllability of Euler equation for rigid body
in \cite{AS86} and was named there {\em reduction} of a {\em
control-affine} system. Here we present a particular version adapted
to our problem. Repeated application of this extension settles
controllability issue  for finite-dimensional Galerkin
approximations of NS systems.

\begin{proposition}\label{lext}
  Consider   control-affine analytic system
\begin{equation}\label{f012}
    \dot{x}=f^0(x)+f^1(x)\hat{v}_1+f^2(x)\hat{v}_2.
\end{equation}
Let
\begin{equation}\label{skob}
 [f^1,f^2]=0, \ [f^1,[f^1,f^0]] =0.
\end{equation}
 Then the system
\begin{equation*}%%\label{f01212}
    \dot{x}=f^0(x)+f^1(x)\tilde{v}_1+f^2(x)\tilde{v}_2+[f^1,[f^2,f^0]]v_{12},
\end{equation*}
is  fixed-time Lie extension of (\ref{f012}). $\Box$
\end{proposition}

{\small\bf Sketch of the proof.} Take Lipschitzian functions $
v_1(t), v_2(t), \ (v_1(0)=v_2(0)=0)$ and change in (\ref{f012})
$\hat{v}_1,\hat{v}_2$ to $\eps^{-1}\dot{v}_1(t)+\tilde{v}_1$ and
$\eps \dot{v}_2(t)+\tilde{v}_2$ respectively. We arrive to the
equation
\begin{equation}\label{f1212}
    \dot{x}=f^0(x)+f^1(x)\tilde{v}_1+f^2(x)\tilde{v}_2+
    \left(\eps^{-1}f^1\dot{v}_1(t)+\eps f^2\dot{v}_2(t)\right).
\end{equation}

Applying the 'reduction formula' from \cite{AS86} or alternatively
'variation of constants' formula of chronological calculus
(\cite{AG78}) one can represent the flow of (\ref{f1212}) as a
composition of the flow $\tilde{P}_t$ of the equation
\begin{equation}\label{vco}
\dot{y}=e^{\ad f^1(y)\eps^{-1}v_1(t)+\ad f^2 \eps v_2(t)}f^0(y)+
f^1(x)\tilde{v}_1+f^2(x)\tilde{v}_2.
\end{equation}
and the flow
\begin{equation}\label{floba}
P_t=e^{f^1\eps^{-1}v_1(t)+f^2 \eps v_2(t)}.
\end{equation}
For the validity of this decomposition the equality $[f^1,f^2]=0$ is
 important.

In (\ref{vco}) $e^{\ad_f}$ is exponential of the operator $\ad_{f}:
\ e^{\ad_f}=\sum_{j=0}^\infty (\ad_f)^j/j!$. The  operator $\ad_{f}$
is determined by the vector field $f$ and acts on vector fields as:
$\ad_f g=[f,g]$, $[f,g]$ being the Lie bracket of $f$ and $g$.

By first of the relations (\ref{skob}) operators $\ad_{f^1}$ and
$\ad_{f^2}$ commute and by the second of these  relations any
iterated Lie bracket of the form $(\ad_{f^{i_1}}) \circ \cdots \circ
(\ad_{f^{i_m}})f^0$, with $i_j=1,2$, vanishes whenever it contains
$\ad_{f^1}$ at least twice.

Taking the expansion of the operator exponential in (\ref{vco}) and
using these facts  we get
\footnote{The time-variant vector field
abbreviated by $O(\eps)$ in (\ref{nper}) equals
$$\eps \phi(\eps\ad_{ v_2(t)f^2}) \ad_{v_2(t)f^2}f^0 +
\eps^2 \ad^2_{v_2(t)f^2} \varphi(\eps\ad_{v_2(t)f^2})[f^1,f^0],
$$ where $\phi(z)=z^{-1}(e^z-1),\ \varphi(z)=z^{-2}(e^z-1-z) $.}
\begin{eqnarray}\label{nper}
\dot{y}=f^0(y)+f^1(x)\tilde{v}_1+f^2(x)\tilde{v}_2+ \\
+\eps^{-1}[f^1,f^0](x)v_1(t)
+[f^1,[f^2,f^0]](x)v_1(t)v_2(t)+O(\eps), \nonumber
\end{eqnarray}

To obtain the flow of (\ref{f1212}) we need to compose the flow of
(\ref{nper}) with the flow (\ref{floba}). For any fixed $T$ one can
get $P_T=Id$ by choosing $v_1(\cdot),v_2(\cdot)$ such that
$v_1(T)=v_2(T)=0$.

From now on we  deal with the fixed $T$ and the flow of the equation
(\ref{nper}).

Let us change  in (\ref{nper})  $v_j(t)$ to
$v_j(t)=2^{1/2}\sin(t/\eps^2)\bar{v}_j(t),$ with $\bar{v}_j(t)$
being functions of bounded variation $j=1,2$. The relaxation
seminorms of the time-variant vector field
$\eps^{-1}[f^1,f^0](x)2^{1/2}\sin(t/\eps^2)\bar{v}_1(t)$ in the
right-hand side of (\ref{nper}) are $O(\eps)$ as $\eps \rightarrow
+0$. Also in the right-hand side of (\ref{nper})
\begin{eqnarray*}
[f^1,[f^2,f^0]](x)2\sin^2(t/\eps^2)\bar{v}_1(t)\bar{v}_2(t)=\\
=[f^1,[f^2,f^0]](x)\bar{v}_1(t)\bar{v}_2(t)-
[f^1,[f^2,f^0]](x)\cos(2t/\eps^2)\bar{v}_1(t)\bar{v}_2(t).
\end{eqnarray*}
The relaxation seminorms of the  addend
$[f^1,[f^2,f^0]](x)\cos(2t/\eps^2)\bar{v}_1(t)\bar{v}_2(t)$ are
$O(\eps^2)$ as $\eps \rightarrow +0$.

Hence the right-hand sides of the equation (\ref{nper}) with the
controls $v_j(t)=2^{1/2}\sin(t/\eps^2)\bar{v}_j(t), \ j=1,2,$
converge in relaxation metric to the vector field
$$f^0(y)+f^1(x)\tilde{v}_1+f^2(x)\tilde{v}_2+  [f^1,[f^2,f^0]](x)\bar{v}_1(t)\bar{v}_2(t)$$
as $\eps \rightarrow 0$. We can consider the product
$\bar{v}_1(t)\bar{v}_2(t) $ as a new control $v_{12}$ and invoke the
Proposition~\ref{fo1}.  $\Box$

\section{Bracket computation in finite and infinite dimensions;
controlling along "principal axes"}\label{pax}
%%%%%%%%%%%%%%%%%%%%%%%%%%%%%%%%%%%%%%%%%%%%%%%%%%%%%%%%%%%%%%%%%
In this section we  adjust the statement of the
Proposition~\ref{fo1} for studying controllability of the systems
(\ref{evort}) and (\ref{galobs}).

From  viewpoint of geometric control the Galerkin approximation
(\ref{galobs}) of NS/Euler system is  particular type of
control-affine system (\ref{caf}). Its state space is
finite-dimensional and is  generated by a finite number of
 eigenfunctions of the Laplace-Beltrami operator $\Delta$ or {\it modes}.

The dynamics of this control system is determined by  quadratic {\it
drift}  vector field
\begin{equation*}%%\label{pdrift}
f_o^0= \Pi^{o}\{\Delta^{-1}q^{o},q^{o}\}+\nu\Delta q^{o}
\end{equation*}
 and by controlled forcing $\sum_{i=1}^r f^i(x)u_i$, where $f^i$ are  constant
($q^{o}$-in\-de\-pen\-dent) cont\-rolled vec\-tor fields.

Start with computation of the particular Lie brackets, which are
involved in the formulation of the Proposition~\ref{fo1}. For two
constant vector fields $f^1,f^2$ there holds
\begin{eqnarray*}
[f^i,f^0_o]= \Pi^{o}\left(\{f^i, \Delta^{-1}w\}+\{w,
\Delta^{-1}f^i\}\right)+ \nu \Delta f^i, \ i=1,2;
\\ %
\! [f^1,[f^2,f^0_o]]=\Pi^{o}\left(\{f^2, \Delta^{-1}f^1\}+\{f^1,
\Delta^{-1}f^2\}\right), \\
\! [f^1,[f^1,f^0_o]]=2\Pi^{o}\{f^1,
\Delta^{-1}f^1\}.
\end{eqnarray*}

This computation is finite-dimensional but the same holds true if
one considers constant vector fields acting in infinite-dimensional
Hilbert space. Taking 'drift' vector field of the equation
(\ref{evort}) in infinite dimensions
\begin{equation*}%%\label{drift}
f^0= \{\Delta^{-1}q,q\}+\nu\Delta q.
\end{equation*}
we conclude with

\begin{lem}\label{f0f1f2} For two constant vector fields $f^1,f^2$ there holds
\begin{eqnarray*}
[f^i,f^0]= \{f^i, \Delta^{-1}w\}+\{w, \Delta^{-1}f^i\}+\nu \Delta
f^i, \ i=1,2;
\\ %
\label{fkk} [f^1,[f^1,f^0]]=2\{f^1, \Delta^{-1}f^1\},
\end{eqnarray*}
\begin{equation}\label{fkl0}
 \mathcal{B}(f^1,f^2)=[f^1,[f^2,f^0]]=\{f^2,
\Delta^{-1}f^1\}+\{f^1, \Delta^{-1}f^2\}. \ \Box%
\end{equation}
\end{lem}

Let us see what is needed for the assumptions of the
Proposition~\ref{lext} to hold. As long as $f^1,f^2$ are constant
and hence commuting, all  we   need is
\begin{equation}\label{pue}
[f^1,[f^1,f^0]]=\{f^1, \Delta^{-1}f^1\}=0.
\end{equation}

 In terms of Euler
equation for {\it ideal} fluid
$$\frac{\partial w}{\partial t}-\{\Delta^{-1}w,w\}=0. $$
this means that $f^1$ corresponds to its {\it steady motion}. In
particular {\it eigenfunctions of the Laplace-Beltrami operator}
$\Delta$ correspond to steady motions and satisfy (\ref{pue}).  In
what follows these eigenfunctions will be used as controlled
directions.

Eigenfunctions of the Laplacian are analogous to principal axes of a
(multidimensional) rigid body.

By Proposition~\ref{lext} given two constant controlled vector
fields one of which corresponds to a steady motion we can extend our
control system by new controlled vector field (\ref{fkl0}) which is
again constant.

Our method consists of iterating this procedure. The
algebraic/geometric difficulties, which arise  on this way, consist
of scrutinizing the newly obtained controlled directions. In
particular we want to know whether at each step we obtain a {\it
proper} extension. This will be illustrated in the further Sections
which deal with particular 2D domains.

Another (analytic) difficulty arises when we pass from
finite-dimensional approximations to controlled PDE. For these
latter the sketched above proof of the Proposition~\ref{lext} is not
valid anymore (for example one can not speak about flows). We have
to reprove statement of the Proposition in each particular
situation. The main idea will be still based on using fast
oscillating control and relaxation metric. The analytic difficulties
are in proving continuity of forcing/trajectory map with respect to
such metric. We will provide a brief comment in the next Section;
the readers can find details in \cite{AS43,AS06,Rod0,Rod1}.

\section{Controllability and accessibility of Galerkin approximations of NS/Euler system on
$\mathbb{T}^2$}
\label{galer}

Here we survey results on accessibility and controllability of
Galerkin approximations of 2D NS/Euler system on $\mathbb{T}^2$.

\subsection{Accessibility of Galerkin approximations}
The result of the computation (\ref{fkl0}) in the periodic case is
easy to visualize when the constant controlled vector fields which
correspond to the eigenfunctions of Laplacian $\Delta$ on
$\mathbb{T}^2$, are written as complex exponentials.

For two different complex eigenfunctions $f^1=e^{ik \cdot x},
f^2=e^{i\ell \cdot x}$ of the Laplacian $\Delta$ on $\mathbb{T}^2$,
where $x \in \mathbb{R}^2, \ k, \ell \in \mathbb{Z}^2$, the Poisson
bracket in (\ref{fkl0}) equals
\begin{equation}\label{245}
\mathcal{B}\left(e^{ik \cdot x}, e^{i\ell \cdot x}\right)=(k \wedge
\ell)(|k|^{-2}-|\ell|^{-2}) e^{i(k+\ell) \cdot x},
\end{equation}
i.e. again corresponds to an eigenfunction of $\Delta$,  provided
that $|k| \neq |\ell|, \ k \wedge \ell \neq 0$. The conclusion is
that given two pairs of complex exponentials $e^{\pm ik \cdot x},
e^{\pm i\ell \cdot x}$ as controlled vector fields we can extend
them by controlled vector fields $e^{i(\pm k \pm \ell) \cdot x}$.

Iterating the computation of the Lie-Poisson  brackets (\ref{fkl0})
and obtaining new directions we end up  with  (finite or infinite)
set of functions, which contains $e^{\pm ik \cdot x}, e^{\pm i\ell
\cdot x}$ and is invariant with
 respect to the bilinear operation $\mathcal{B}(\cdot , \cdot)$.

Therefore in the case of $\mathbb{T}^2$ starting with   controlled
vector fields, which
 correspond to the eigenfunctions  $e^{ik \cdot x}, \ k \in \mathbb{Z}^2,$ of the Laplacian, the whole
 computation
 of Lie extensions 'can be modeled' on the integer lattice $\mathbb{Z}^2$
 of "mode indices" $k$.

Actually one has to operate with real eigenfunctions of the
Laplacian on $\mathbb{T}^2$, i.e. with the functions of the form
$\cos (k \cdot x), \sin (k \cdot x)$. Also in this case computation
of iterated Lie-Poisson brackets (\ref{fkl0}) can be modelled on
$\mathbb{Z}^2$ and the addition formulae are similar to those for
the complex case.

We formulate now bracket generating criterion.

\begin{proposition}[bracket generating property]
\label{bragen}
 If
\begin{equation}\label{kl}
 |k| \neq |\ell|, \ | k \wedge \ell |=1,
\end{equation}
  then:

i)  invariant with
 respect to $\mathcal{B}$ set of functions, which contains complex functions $e^{\pm ik \cdot x}, e^{\pm i\ell \cdot
 x},$
 contains all the eigenfunctions $e^{im \cdot x}, \ m \in \mathbb{Z}^2 \setminus
 0$;

ii) invariant with
 respect to $\mathcal{B}$ set of real functions, which contains $$\cos (k
\cdot x),  \sin (k \cdot x), \cos (\ell  \cdot x),\sin (\ell \cdot
x),$$ contains all the eigenfunctions $\cos (m  \cdot x),\sin (m
\cdot x), \  m \in \mathbb{Z}^2 \setminus
 0. \ \Box$
\end{proposition}

Bracket generating property for Galerkin approximation of 2D NS
system under periodic boundary conditions has been established  by
W.E and J.Mathingly in \cite{EM}.
 The following result from \cite{EM} is an immediate consequence of
the previous Proposition~\ref{bragen} and Theorem~\ref{kren}.

\begin{corollary}[accessibility by means of 4 controls]
\label{ac4gal} For any subset $\mathcal{M} \subset \mathbb{Z}^2$
there exists a larger set $\mathcal{M}' \supseteq \mathcal{M}$ such
that Galerkin $\mathcal{M}'$-approximation controlled by the forcing
\begin{equation}\label{acc4}
\cos(k \cdot x)v_k(t)+\sin(k \cdot x)w_k(t)+ \cos(\ell \cdot x)
v_\ell(t)+\sin(\ell \cdot x) w_\ell(t),
\end{equation}
with $k,\ell$ satisfying (\ref{kl}), is strongly accessible. $\Box$
\end{corollary}

Here $4$ controls $v_k(t),w_k(t),v_\ell(t),w_\ell(t) $ are used  for
providing strong accessibility, but actually it can be achieved by a
smaller number of controls.

{\small\bf Example (accessibility by means of 2 controls)} Consider
 forcing
\begin{equation}\label{acc2}
g v(t)+\bar{g} \bar{v}(t), \ g=\cos(k \cdot x)+\cos(\ell \cdot x), \
\bar{g}=\sin(k \cdot x)-\sin(\ell \cdot x).
\end{equation}
 Application of controlled forcing
(\ref{acc4}) results in $4$ independent controls, each one
$v_k(t),w_k(t),v_\ell(t),w_\ell(t) $ entering one  and only one of
the equations (\ref{coord}). Under the action of (\ref{acc2}) each
of controls $v,\bar{v}$ is applied simultaneously to a pair of
equations (\ref{coord}).

Assume $|k| \neq |\ell|, \ k \wedge \ell \neq 0$. We compute the
bilinear form (\ref{fkl0}).
\begin{equation*}
\mathcal{B}(g, g) =(-|\ell|^{-2}+|k|^{-2})\{\cos(k \cdot
x),\cos(\ell \cdot
 x)\}
\end{equation*}
Up to a scalar multiplier $\mathcal{B}(g, g)$ equals
\begin{eqnarray*}
 (-|\ell|^{-2}+|k|^{-2})\sin(k \cdot x)\sin(\ell
\cdot
 x)= \nonumber\\ =( k \wedge \ell ) (-|\ell|^{-2}+|k|^{-2})\left(\cos((k-\ell) \cdot
 x)- \cos((k+\ell) \cdot
 x)\right).
 \label{gg}
\end{eqnarray*}
Similarly up to a scalar multiplier $\mathcal{B}(\bar{g}, \bar{g})$
equals
\begin{equation*}
  ( k \wedge \ell ) (-|\ell|^{-2}+|k|^{-2})\left(\cos((k-\ell) \cdot
 x)+ \cos((k+\ell) \cdot
 x)\right);
\end{equation*}
The span of $\mathcal{B}(g, g), \mathcal{B}(\bar{g}, \bar{g})$
coincides with the span of
$$g^{01}=\cos((k-\ell) \cdot
 x), \  g^{21}=\cos((k+\ell) \cdot
 x).$$

The direction $$\bar{g}^{01}=\sin((k-\ell) \cdot  x)$$ is obtained
from computation of $\mathcal{B}(g, \bar{g})$.

Choose $k=(1,1), \ell=(1,0)$; then  $m=k+\ell=(2,1),
n=k-\ell=(0,1)$.

Computing new directions $\mathcal{B}(g^{01}, \bar{g}),
\mathcal{B}(\bar{g}^{01},g )$ we
 note  that due to the equality $|n|=|\ell|$ they coincide with
$\mathcal{B}(g^{01},\sin(k \cdot x) ),
\mathcal{B}(\bar{g}^{01},\cos(k \cdot x) ) $ respectively, and their
span coincide with the span of functions
$$\bar{g}^{12}=\sin((k+n) \cdot
x) ), \bar{g}^{10}=\sin((k-n) \cdot x) )= \sin (\ell \cdot x) )$$

Similarly span of the directions
$$\mathcal{B}(g^{01}, g), \mathcal{B}(\bar{g}^{01},\bar{g} )$$
coincides with the span of functions
$$g^{12}=\cos((k+n) \cdot
x) ), g^{10}=\cos( (k-n) \cdot x)=\cos(\ell \cdot x)  .$$

Then $g-g^{10}=\cos(k \cdot x), \ \bar{g}-\bar{g}^{10}=\sin(k \cdot
x)$.

These two  functions together with the functions $g^{01},
\bar{g}^{01}$ form a quadruple, which satisfies the assumptions of
the Corollary~\ref{ac4gal}. Hence our system is accessible by means
of $2$ controls.  $\Box$

\begin{remark} It is plausible that strong accessibility of Galerkin approximation can be
achieved by one control. $\Box$
\end{remark}

\subsection{Controllability of Galerkin approximations}

What regards controllability, then in general bracket generating
property is not sufficient for it. One has to select Lie brackets,
which lead to Lie extensions, meanwhile in the previous example
$\{g, \Delta^{-1}g\},\{\bar{g}, \Delta^{-1}\bar{g}\}$ a priori do
not correspond to Lie extensions.

Even in the finite-dimensional case one needs a stronger result of
Proposition~\ref{lext} to prove controllability property for the
Galerkin approximations. This has been accomplished in
\cite{ASpb,Rom} in @D and 3D cases.

\begin{thm}\label{acgal}
Let $k, \ell$ satisfy (\ref{kl}). For any subset $\mathcal{M}
\subset \mathbb{Z}^2$ there exists a larger set $\mathcal{M}'$ such
that Galerkin $\mathcal{M}'$-approximation controlled by the forcing
\begin{equation*}
\cos(k \cdot x)v_k(t)+\sin(k \cdot x)w_k(t)+ \cos(\ell \cdot x)
v_\ell(t)+\sin(\ell \cdot x) w_\ell(t),
\end{equation*}
 is globally controllable. $\Box$
\end{thm}

The proof of this controllability result consists of iterated
application of Lie extension described in the
Proposition~\ref{lext}.  At each step we extend the system by new
controlled vector fields  which according to {fkl0}) and
(\ref{245}) correspond to $f^{m \pm \ell}=\cos((m \pm \ell) \cdot
x), \ \bar{f}^{m \pm \ell}=\sin((m \pm \ell) \cdot x)$.

At the end of the iterated procedure we arrive to the  system with
extended set of controls - one for each observed mode. This  latter
system evidently would be controllable.

An important case when controllability  of Galerkin approximation is
implied by bracket generating property regards  2D Euler equation
for incompressible ideal fluid ($\nu=0$).

Indeed in this case the drift (zero control) dynamics is {\it
Hamiltonian} and it evolves on a compact energy level. By Liouville
and Poincare theorems Poisson-stable points of this dynamics are
dense and one can apply Lobry-Bonnard theorem (\cite{ASkv,Ju}to
establish

\begin{thm}\label{bonn}
For $\nu=0$ the Galerkin approximation of the 2D Euler system
controlled is globally controllable by means of the forcing
(\ref{acc2}). $\Box$
\end{thm}

We trust that in the case of ideal fluid the controllability of
Galerkin approximation of 2D Euler system can be achieved by {\it
scalar} control.

\section{Steady state controlled directions: abstract controllability result for NS system}
\label{steady}

We can not apply Proposition~\ref{lext} to infinite dimensions
directly. Still the main idea of adding new controlled directions is
valid for NS systems. Now we want to formulate an abstract
controllability criterion, based on Lie extensions and on the
computation of Lemma~\ref{f0f1f2}. This criterion will be employed
in forthcoming Sections for establishing controllability of NS
system on various  2D domains.

\begin{thm}[controllability of NS systen via saturation of controls]
\label{61} Let $$\spn\{f^1, \ldots ,
f^r\}=\mathcal{S}=\mathcal{D}^1$$ be a finite dimensional space of
controlled directions. Assume $f^1, \ldots f^r$ to be steady motions
of Euler equation (\ref{pue}).  For each pair of linear subspaces
$\mathcal{L}^1, \mathcal{L}^2$ consider the span of the image
$\mathcal{B}(\mathcal{L}^1, \mathcal{L}^2)$ of the bilinear map
(\ref{fkl0}). Define successively
\begin{equation*}%%\label{dit}
\mathcal{D}^{j+1}=\mathcal{D}^{j}+ \spn \mathcal{B}(\mathcal{S}^{j},
\mathcal{D}^{j}), \ j=1,2, \ldots
\end{equation*}
where $\mathcal{S}^{j} \subseteq  \mathcal{D}^{j}$ is a linear
subspace spanned by steady motions.  If  $\bigcup_j \mathcal{D}^{j}$
is dense in the Sobolev space  $H^2(M)$, then the NS system is
controllable in finite-dimensional projections and
$L_2$-approximately controllable. $\Box$
\end{thm}

As far as $\mathcal{D}^{1}$ consists of steady motions then
$\mathcal{D}^{j+1} \supseteq \mathcal{D}^{j}+ \spn
\mathcal{B}(\mathcal{D}^{1}, \mathcal{D}^{j}) .$ Introduce the
sequence of spaces
\begin{equation}\label{dlit}
    \mathcal{D}_1^{j+1}=\mathcal{D}_1^{j}+ \spn
\mathcal{B}(\mathcal{D}^{1}, \mathcal{D}_1^{j}).
\end{equation}
Evidently $\mathcal{D}_1^{j} \subseteq \mathcal{D}^{j}$ and density
of $\bigcup_j \mathcal{D}_1^{j}$ in $H^2(M)$ guarantees
controllability.

Let for $f_s \in \mathcal{D}^{1}, \ D_{f_s}=\mathcal{B}(f_s,
\cdot)$:
$$D_{f_s}=\{\Delta^{-1}\cdot,f_s\}+\{\Delta^{-1}f_s,\cdot\}.$$
The iterated computations (\ref{dlit}) correspond to iterated
applications of the operators $D_{f_s}$ to $f^1, \ldots f^r$ and
taking linear span.
\begin{corollary}\label{fdens}   Let $\mathcal F$ be the minimal common invariant
linear subspace of the operators $D_{f_1},\ldots,D_{f_l}$ which
contains $f_1,\ldots,f_k$. If $\mathcal F$  is everywhere dense in
$L_2(M)$, then the system is $L_2$-appro\-xi\-ma\-te\-ly
controllable and solidly controllable in finite dimensional
projections. $\Box$
\end{corollary}

\section{Navier-Stokes and Euler System on $\mathbb{T}^2$}\label{t2}

In this Section we formulate results regarding  controllability in
finite-dimensional projections and $L_2$-approximate controllability
on $\mathbb{T}^2$. Namely we describe sets of controlled directions
which satisfy criterion provided by Theorem~\ref{61}.

We take  the basis of complex eigenfunctions $\left(e^{ik \cdot
x}\right)$ of the Laplacian on $\mathbb{T}^2$ and take the Fourier
expansion of the vorticity $w(t,x)=\sum_k q_k(t)e^{ik \cdot x}$ and
control $v(t,x)=\sum_{k \in \mathcal{K}^1} v_k(t)e^{ik \cdot x}$;
here $k \in \mathbb{Z}^2$. As far as $w$ and $f$ are real-valued, we
have $\bar{q}_n=q_{-n}, \ \bar{v}_n=v_{-n}$. We assume $v_0=0$ and
$q_0=0$.

Using the computation (\ref{245}) we write the infinite system of
ODE (\ref{coord}) as
\begin{equation}\label{mleqn}
    \dot{q}_k=\sum_{m+n=k, |m|<|n|}(m \wedge n)(|m|^{-2}-|n|^{-2})q_m
q_n- |k|^2q_k+\hat{v}_k(t).
\end{equation}

The controls $v_k$ are nonvanishing  only in the equations for the
variables $q_k$ indexed by symmetric set $\mathcal{K}^1 \subset
\mathbb{Z}^2 \setminus \{0\}$. For $k \not\in \mathcal{K}^1$ the
dynamics is
\begin{equation}\label{unc}
\dot{q}_k=\sum_{m+n=k, |m|<|n|}(m \wedge n)(|m|^{-2}-|n|^{-2})q_m
q_n-|k|^2q_k, \ k \not\in \mathcal{K}^1.
\end{equation}

There is a symmetric set of observed modes $\mathcal{K}^{o} \supset
\mathcal{K}^1$, which we want to steer to some preassigned values.
In the only interesting case where $\mathcal{K}^1$ is proper subset
of $\mathcal{K}^{o}$, the equations indexed by $k \in
\mathcal{K}^{o} \setminus \mathcal{K}^{1}$ are of the form
(\ref{unc}). They do not contain controls and need to be {\it
controlled via state variables}.

We give a hint of how this can be done; it is an
infinite-dimensional version of Proposition~\ref{lext} for NS system
on $\mathbb{T}^2$.

 Let $r,s \in \mathcal{K}^1, \ r \wedge s \neq 0, |r| <
|s|$, and $k=r+s \not\in \mathcal{K}^1$. The equations for $q_r,q_s$
contains controls $\hat{v}_r,\hat{v}_s$, while the equation for
$q_k$ does not.

Take Lipschitzian functions $ v_r(t), v_s(t), \ (v_r(0)=v_s(0)=0)$
and substitute $\hat{v}_r,\hat{v}_s$ in the equations for $q_r,q_s$
 by $\eps^{-1}\dot{v}_r(t)+\tilde{v}_r$ and $\eps
\dot{v}_s(t)+\tilde{v}_s$. We obtain the equations
\begin{eqnarray*}
\dot{q}_r=\sum_{m+n=r, |m|<|n|}(m \wedge n)(|m|^{-2}-|n|^{-2})q_m
q_n- |r|^2q_r+\eps^{-1}\dot{v}_r(t)+\tilde{v}_r;
 \\
\dot{q}_s=\sum_{m+n=s, |m|<|n|}(m \wedge n)(|m|^{-2}-|n|^{-2})q_m
q_n- |s|^2q_s+\eps \dot{v}_s(t)+\tilde{v}_s.
\end{eqnarray*}

Introduce variables $q^*_r=q_r-\eps^{-1}v_r(t),q^*_s=q_s-\eps v_s(t)
$. Assuming $v_r(T)=v_s(T)=0$ we conclude
$q_r(T)=q^*_r(T),q_s(T)=q^*_s(T)$.

Let us rewrite the infinite system of ODE (\ref{mleqn})-(\ref{unc})
via $q^*_r,q^*_s$ in place of $q_r,q_s$. The right-hand side of the
equation for $q_k=q_{r+s}$ gains the addend
$$(r \wedge s)(|r|^{-2}-|s|^{-2})(q^*_r+\eps^{-1}v_r(t))(q^*_s+\eps
v_s(t))$$ and we see that the  controls $v_r, v_s$ enter this
equation via product $v_r(t)v_s(t)$. Same $v_r,v_s$  enter this and
all other equations linearly.

Substitute   $v_j(t)$ for $j=r,s$ by
$v_j(t)=2^{1/2}\sin(t/\eps^2)\bar{v}_j(t)$ with $\bar{v}_j(t)$
having bounded variations. Then the right-hand side of the equation
 for $q_k$ will gain the product
$$2\sin^2(t/\eps^2)\bar{v}_r(t)\bar{v}_s(t)=(1
-2\cos(2t/\eps^2))\bar{v}_r(t)\bar{v}_s(t).$$

If $\eps \rightarrow +0$ this product tends to
$\bar{v}_r(t)\bar{v}_s(t)$ in relaxation metric. In all other
equations $\bar{v}_r(t)$ and  $\bar{v}_s(t)$ enter linearly and are
multiplied by fast oscillating functions. Therefore the
corresponding terms tend to $0$ in relaxation metric.

Therefore one can pass (as $\eps \rightarrow 0$) to a limit system
which now contains 'new' control
$\bar{v}_{rs}=\bar{v}_r(t)\bar{v}_s(t)$ in the equation for
$q_k=q_{r+s}$. (This control corresponds to the control $v_{12}$
from the Proposition \ref{lext}).

A difficult analytic part  is justification of this passage to the
limit. It is accomplished in \cite{AS43,AS06} for $\mathbb{T}^2$ and
in \cite{Rod0,Rod1} for  rectangular and other kinds of regular 2D
domains. We refer  interested readers to these publications.

Note that the new controlled direction corresponds to complex
exponential which is an eigenfunction of the Laplacian on
$\mathbb{T}^2$. Hence we can model Lie extensions and formulate
controllability results in terms of indices $k \in \mathbb{Z}^2$ of
controlled modes.

Define a sequence of sets $\mathcal{K}^j \subset \mathbb{Z}^2$
iteratively as follows:
\begin{eqnarray}\label{setk2}
j=2, \ldots  , \\
\mathcal{K}^j  =\mathcal{K}^{j-1} \bigcup   \{m+n | \ m,n \in
\mathcal{K}^{j-1} \bigwedge \|m\| \neq \|n\| \bigwedge m \wedge n
\neq 0 \}. \nonumber
\end{eqnarray}

\begin{definition}
\label{sats} A finite set $\mathcal{K}^1 \subset \mathbb{Z}^2
\setminus \{0\}$ of forcing modes is called saturating if
$\bigcup_{j=1}^\infty \mathcal{K}^j=\mathbb{Z}^2 \setminus \{0\}$,
 where $\mathcal{K}^j$ are defined by (\ref{setk2}). $\Box$
\end{definition}

\begin{thm}{\it (controllability in finite-dimensional
projection)} \label{procon}
 Let $\mathcal{K}^1$ be a saturating  set of controlled forcing
 modes and $\mathcal{L}$ be any finite-dimensional subspace of $H^2(\mathbb{T}^2)$.
  Then for any  $T>0$ the  NS/Euler system on $\mathbb{T}^2$
is time-$T$ solidly controllable in finite-dimensional projections
and time-$T$ $L_2$-approximately controllable. $\Box$
\end{thm}

As we see the saturating property is crucial for controllability. In
\cite{AS06} the following  characterization of this property has
been established.

\begin{thm}\label{snsat}
For a symmetric finite set $\mathcal{K}^1=\{m^1, \ldots , m^s\}
\subset \mathbb{Z}^2$ the following properties are equivalent:

i) $\mathcal{K}^1$ is saturating;

ii) the greatest common divisor of the numbers $d_{ij}=m^i \wedge
m^j, \ i,j \in \{1,\ldots , s\}$ equals $1$  and  there exist
$m^\alpha, m^\beta \in \mathcal{K}^1$, which are not collinear and
have different lengths. $\Box$
\end{thm}

\begin{corollary}
The set $\mathcal{K}^1=\{(1,0),(-1,0),(1,1),(-1,-1)\} \subset
\mathbb{Z}^2$ is saturating. Solid controllability in any
finite-dimensional projection and $L_2$-approximately
controllability can be achieved by forcing $4$ modes. $\Box$
\end{corollary}

\section{Controllability of 2D NS system on a rectangular domain}
\label{rect}

The study of controllability in finite-dimensional projections and
of $L_2$-approximate controllability on a rectangular domain has
been accomplished  by S.Rodrigues in \cite{Rod0,Rod1}. The main idea
is similar to the one employed for the periodic case, but the
computations are more intricate. The reason is twofold: i) the
algebraic properties of the bilinear operation  calculated for the
eigenfunctions of the Laplacian are more complex, and ii) one needs
to care about boundary conditions.

For a velocity field $u$ on a rectangular $\mathcal{R}$ with sides
of lengths $a,b \ (a \neq b)$ we assume Lions boundary conditions to
hold. In terms of the vorticity $w$ they can be written as
(\ref{bc}).

The (vorticity) eigenfunctions $\phi^k$ of the Laplacian are
\begin{equation}\label{recteig}
\phi^k=\sin \left(\frac{\pi}{a}k_1x_1\right)\sin
\left(\frac{\pi}{b}k_2x_2\right), \ (k=(k_1,k_2) \in \mathbb{Z}^2.
\end{equation}

To find  extending controlled direction one needs to pick two
eigenfunctions $f^1=\phi^k, f^2=\phi^\ell, \ k,\ell \in
\mathbb{Z}^2$ and to proceed with the computation (\ref{fkl0}). This
results in a linear combination of at most $4$ eigenfunctions $W^s$.

Then again one can follow  Lie extensions on two-dimensional lattice
$\mathbb{Z}^2$  of  Fourier exponents $k=(k_1,k_2)$. If the
controlled modes are indexed by $ k \in \mathcal{K}^1=\{(k_1,k_2)| \
1 \leq k_1, k_2 \leq 3, \ k\neq (3,3)\}$, then one can verify that
after $m$ Lie extensions the set of extended controlled directions
will contain all the modes $(k_1,k_2)$ with $k_1,k_2 \leq m+3$ with
the exception of $(m+3,m+3)$.

This leads to  the following controllability result.
\begin{thm}[controllability on rectangular
domain]\label{rectan}  Let $8$  controlled directions correspond to
the functions (\ref{recteig}) with $ k \in \{(k_1,k_2)| \ 1 \leq
k_1, k_2 \leq 3, \ k\neq (3,3)\}$.Then the NS system defined on the
rectangular domain under Lions boundary condition is controllable in
finite-dimensional projections and $L_2$-approximately controllable.
$\Box$
\end{thm}

\begin{remark}
In the previous Theorem the controllability in projections and
approximate controllability are guaranteed by $8$ controlled
directions. $\Box$
\end{remark}

\section{Controllability in finite-dimensional projection of 2D NS/Euler system on a smoothened
rectangular domain}\label{smooth}

We are going to prove that the established above controllability of
NS system on a rectangular domain persists if one slightly perturbs
the rectangular and in particular smoothens it.

To simplify the presentation we assume the rectangular $\mathcal{R}
$ to have dimensions $a,b$ with $a/b$  irrational. The irrationality
of $a/b$ guarantees that the eigenvalues of the Laplacian on
$\mathcal{R} $ are simple. Let $A(0;0)$, $B(0;b)$, $C(a;b)$,
$D(0;a)$ be the vertices of the rectangular.

Consider parameterized quadruples of points
$$A_\eps(\eps;\eps), B(\eps;b-\eps), C(a-\eps;b-\eps),
D(\eps;a-\eps).$$ Construct $\eps$-circumferences centered at these
points and round off the corners of $\mathcal{R} $ by the quarters
of these circumferences. We obtain smoothened rectangular
$\mathcal{R}_\eps$.

The following results on continuous dependence of eigenvalues and
eigenfunctions of the Laplacian on the domain holds (see e.g.
\cite[Vol. I,Ch.VI]{CH}): for each $N$ choosing $\eps >0$
sufficiently small one can maintain first $N$ eigenvalues
$\lambda_j^\eps$ of $\Delta$ on $\mathcal{R}_\eps$ simple; besides
the eigenspace which corresponds to $j$-th eigenvalue
$\lambda_j^\eps$ of $\Delta$ on $\mathcal{R}_\eps$ converges to the
respective eigenspace of $\Delta$ on $\mathcal{R}$ as well as
$\lambda_j^\eps \rightarrow \lambda_j$, as $\eps \rightarrow 0.$

Some more technical problems arise from the fact that due to
boundary conditions the evolution space of NS system depends on the
boundary and one has to take this dependence into account.

\begin{thm}[solid controllability on smoothened rectangular
domain] Consider NS system on rectangular domain $\mathcal{R}=[0,a]
\times [0,b]$ ($b/a$ - irrational) controlled along directions which
correspond to the eigenfunctions
\begin{equation}\label{36}
 \sin \left(\frac{\pi}{a}k_1x_1\right)\sin
\left(\frac{\pi}{b}k_2x_2\right), | \  \leq k_1, k_2 \leq 3, \ k\neq
(3,3)\}.
\end{equation}
 Let $\mathcal{R}_\eps$
be the described above perturbations of the rectangle $\mathcal{R}$.
Choose for each  $\eps
>0$ eigenfunctions of the Laplacian on $\mathcal{R}_\eps$  (and
respective controlled directions which are perturbations of the
functions (\ref{36}). Then for any finite-dimensional space
$\mathcal{L}$ there exists $\eps_\mathcal{L}>0$ such that  $\forall
\eps \in (0, \eps_\mathcal{L}]$ the NS system on the perturbed
rectangle $\mathcal{R}_\eps$  controlled by these directions is
controllable in projection on $\mathcal{L}. \ \Box$
\end{thm}

\begin{remark}
Further on  there is a need to deal with  domains with analytic
boundary. It is evident that one can slightly perturb the smoothened
rectangles $\mathcal{R}_\eps$  to obtain domains with analytic
boundaries (also denoted by $\mathcal{R}_\eps$)  for which the claim
of the previous Theorem still holds. $\Box$
\end{remark}

\section{Controllability on a
generic 2D domain homeomorphic to a disc}\label{generdom}

In this section we consider NS system on a 2D domain $M$ under
boundary conditions (\ref{bc}). We manage to prove that for a
generic domain (exact meaning of genericity will be specified in a
moment) one can
 choose $8$ controlled directions (corresponding to eigenfunctions
 of the Laplacian $\Delta$ on $M$) which guarantee
controllability in finite-dimensional projections and
$L_2$-approximate controllability.

In what follows we assume $M$ to have analytic boundary, to be
 homeomorphic to a disc and to be endowed with Euclidean metric.
 Given another such domain $D$ one can map $M$ on $D$ analytically
 (Riemann mapping theorem). The induced Riemannian metric on $D$
can be represented (see \cite[Vol. 1, \S\S 11-13]{DNF} for an
elementary account) in conformal form as:
\begin{equation}\label{rimet}
\mu=e^{a(x_1,x_2)}(dx_1^2+dx_2^2).
\end{equation}

Its curvature computed as $K=(1/2)e^{-a(x_1,x_2)}(\partial^2
a/\partial x_1^2+\partial^2 a/\partial x_2^2)$ must vanish, i.e.
\begin{equation}\label{plos}
\partial^2 a/\partial x_1^2+\partial^2 a/\partial x_2^2=0.
\end{equation}
On the contrary if $D$ possesses Riemannian metric $\mu$ of form
(\ref{rimet}), which satisfies (\ref{plos}), then $D$ can be
isometrically and analytically mapped onto a 2D domain  $M'$ with
Euclidean metric.

Therefore instead of speaking of various 2D domains we will speak
about various metrics (\ref{rimet}) satisfying (\ref{plos}) on a
given domain.

From now on $D$ will be chosen from set of  analytically smoothened
rectangulars for which controllability has been established in the
previous section. Let $\mathcal{R}_0=\mathcal{R}_{\eps_0}$; assume
that $\mathcal{R}_0$ is endowed with a metric (\ref{rimet}), for
which (\ref{plos}) holds. Generic domain in $\mathbb{R}^2$
corresponds to a generic metric (\ref{rimet})  on $\mathcal{R}_0$ or
equivalently to a generic function $a(x_1,x_2)$ satisfying
(\ref{plos}).

Recall that for the metric (\ref{rimet}) the Laplacian $\Delta u$ is
calculated as
\begin{equation}\label{alap}
\Delta u=e^{-a(x_1,x_2)}(\partial^2 u/\partial x_1^2+\partial^2
u/\partial x_2^2).
\end{equation}

Assume that $f_1,\ldots,f_l$ are eigenfunctions of this Laplacian:
$\Delta^{-1}f_s=\lambda^{-1}_s f_s,\quad s=1,\ldots,l. $

Abstract controllability criterion formulated as
Corollary~\ref{fdens} claims that  it suffices to verify that the
iterated applications of
$D_{f_s}=\{\Delta^{-1}\cdot,f_s\}+\{\Delta^{-1}f_s,\cdot\}$ to $f_j$
result in a dense subset of $H^2(M)$.

\begin{thm}\label{102}
For a residual subset of Riemannian discs $M$ with analytic
boundaries there exist $8$ eigenfunctions (modes) $f_1,\ldots,f_8$
of the Laplace-Beltrami operator $\Delta$ on $M$ such that the NS
system on $M$ is controllable in finite-dimensional projections and
$L_2$-approximately controllable by means of (controlled) forcing
applied to the modes $f_1,\ldots,f_8$. $\Box$
\end{thm}

We call a subset of topological space $\mathcal{T}$ {\it residual}
if it contains  an intersection of countable family of open dense
subsets of $\mathcal{T}$.

{\it Sketch of the proof:}

For establishing the two types of controllability it suffices to
establish controllability in projection on any finite-dimensional
{\it coordinate} subspace $\mathcal{L}$ (see \cite{AS06,AS43}.
Controllability in projection on fixed coordinate subspace
$\mathcal{L}$ amounts to  verification of nonnullity of some
determinant, denoted $\DL$, which is calculated  via eigenfunctions
$f_1,\ldots,f_8$ of $\Delta$  and via their iterated Poisson
brackets (Corollary~\ref{fdens}).

For the domains $\mathcal{R}_\eps$ constructed in the previous
Section and endowed with Euclidean metric this determinant $\DL$ is
nonvanishing, provided that $\eps >0$ is sufficiently small (say
$\eps \in (0, \eps_0]$).

Taking any domain $M$ with analytic boundary and Euclidean metric we
establish its isometry to $\mathcal{R}_{\eps_0}$ provided with
metric $\mu_1=e^{a(x_1,x_2)}(dx_1^2+dx_2^2)$ which satisfies
(\ref{plos}).

Consider a 'homotopy'
$$\mu_t=e^{ta(x_1,x_2)}(dx_1^2+dx_2^2), \ 0 \le t \le 1.$$
All the metrics $\mu_t$ are locally euclidean, given that
$\partial^2 a/\partial x_1^2+\partial^2 a/\partial x_2^2=0$.

Obviously the dependence on $t$ of the corresponding Laplacian
$$\Delta(t)=e^{-ta(x_1,x_2)}(\partial^2 /\partial x_1^2+\partial^2
/\partial x_2^2)$$  is analytic.

We want to trace the evolution of a finite number of its eigenvalues
$\lambda_j^t$  $(j \in J \mbox{- finite set})$ and of the
corresponding eigenfunctions with $t$ varying in $[0,1]$.

By classical result of perturbation theory   (see \cite[Ch. 2, Ch.
7]{KaP}) eigenvalues $\lambda^t_j$ of an analytic family of linear
operators are analytic with respect to $t$ beyond finite number of
exceptional points in $[0,1]$. Any moment $t$ for which all the
eigenvalues $\lambda_j^t, \ (j \in J)$  are distinct is
nonexceptional. For generic operator singularities of the function
$t \mapsto \lambda_j^t$  may occur when $\lambda_j^t$  become
multiple. The eigenvectors and respective eigenprojections may have
poles at the exceptional points.

The picture is much more regular for normal operators and in
particular for Laplacians, which are self-adjoint. In this case the
eigenvalues and eigenfunctions are known to depend analytically on
$t$ everywhere on $[0,1]$ (\cite[Ch.2, Th.1.10]{KaP}). Also the
dependence of the derivatives of the eigenfunctions on $t \in [0,1]$
is analytic. Hence the determinant $\DL$ is analytic function of
$t$. As far as it is nonvanishing for $t=0$, it may vanish only at
finite number of points $t \in [0,1]$.

Take $\mu_t$ corresponding to all nonexceptional $t \in [0,1]$ for
which $\DL$ is nonvanishing. Among those there exist $t_s$
arbitrarily close to $1$ and metrics $\mu_{t_s}$ arbitrarily close
to $\mu_1$ for which the eigenvalues of interest are distinct and
$\DL$ is nonvanishing. The dependence of the eigenfunctions and of
their derivatives on metric $\mu$ is then continuous in a
neighborhood of $\mu_{t_s}$. Hence $\DL$ is nonzero for all $\mu$
from small neighborhoods of these $\mu_{t_s}$. Taking union of these
neighborhoods we get an open set whose closure contains $\mu_1$.
Repeating the homotopy argument  for each $\mu_1$ running over the
set of all metrics \ref{rimet})-(\ref{plos}) on
$\mathcal{R}_{\eps_0}$ we get an open dense set of metrics for which
$\DL$ is nonvanishing.

Finally running over the (countable family of all) coordinate
subspaces $\mathcal L$ in $H^2(M)$ and taking the intersection of
the corresponding open sets of metrics we obtain a residual set of
domains for which all the determinants $DL$ are nonvanishing and
hence the assumption of Corollary~\ref{fdens} is verified.
\begin{remark}
The construction of a residual set of Riemannian metrics can be
transferred (almost) without alterations to torus $\mathbb{T}^2$,
for which we studied controllability of NS/Euler system in the
Section~\ref{t2}. The conclusion would be that there exists a
residual set of Riemannian metrics on $\mathbb{T}^2$ such that the
assumptions of the Corollary~\ref{fdens} are verified and therefore
the NS system is controllable in finite-dimensional projections and
$L_2$-approximately controllable on $\mathbb{T}^2$ endowed with any
of these metrics. $\Box$
\end{remark}

\section{NS/Euler system on sphere $\mathbb{S}^2$}
\label{sphere}

The controlled vector fields we employ in the case of $\mathbb{S}^2$
correspond to the eigenfunctions of the corresponding {\it spherical
Laplacian} or to so called {\it spherical harmonics}. We start with
their brief description.

\subsection{Spherical Harmonics}\label{sfegam}
In this subsection  we introduce some notions and results regarding
spherical harmonics; our source was mainly the book \cite[Ch.
10,11]{Arn2} by V.I.Arnold.

Consider sphere $\mathbb{S}^2$ equipped  with the Riemannian metrics
inherited from $\mathbb{R}^3$ and with  area 2-form $\sigma$. the
latter defines symplectic structure on $\mathbb{S}^2$.

The eigenfunctions of the spherical Laplacian are described by the
following classical result. Recall that a function $g$ is
homogeneous of degree $s$ on $\mathbb{R}^n \setminus 0 $, if
$g(\kappa x)=\kappa^s g(x)$ for each $\kappa >0$. A function $g$ is
harmonic in $\mathbb{R}^n \setminus 0 $ if $\Delta g=0$; where
$\Delta$ is the euclidean Laplacian. It is known that a harmonic
homogeneous function of degree $s>0$ is extendable by continuity
($g(0)=0$) to a harmonic function on  $\mathbb{R}^n$. This harmonic
function is smooth and therefore it must be homogeneous polynomial
of integer degree $s>0$.

\begin{thm}[\cite{Arn2}] Constants are eigenfunctions of spherical
Laplacian (of degree $0$). If a (smooth) harmonic function defined
on $\mathbb{R}^n \setminus 0 $ is homogeneous of degree $s>0$, then
its restriction onto sphere is eigenfunction of the spherical
Laplacian $\tilde{\Delta}$ with the eigenvalue $-s(s+n-2)$. Vice
versa every eigenfunction of $\tilde{\Delta}$ is a restriction onto
$\mathbb{S}^n$ of a homogeneous polynomial. $\Box$
\end{thm}

Another famous result is the Maxwell's theorem (\cite{Arn2}), which
holds in $\mathbb{R}^3$. It states that if
$\rho(x)=(x_1^2+x_2^2+x_3^2)^{-1/2}$ is a fundamental solution of
the Laplace equation in $\mathbb R^3$, then any spherical harmonic
$a$ on $\mathbb{S}^2$ can be represented as iterated directional
derivative of $\rho$:
$$
a=l_1\circ\cdots\circ l_n\rho,
$$
where $l_1,\ldots,l_n\in\mathbb R^3$ and the set
$\{l_1,\ldots,l_n\}$ is uniquely determined by $a$.

Our controlled directions will correspond to spherical harmonics on
$S^2$, which are the restrictions to $\mathbb{S}^2$ of homogeneous
functions on $\mathbb R^3$. In particular we invoke so called zonal
spherical harmonics, which are iterated directional derivatives of
$\rho$ with respect to a fixed direction $l$.

Let $a,b$ be smooth (not necessarily homogeneous) functions on
$\mathbb R^3$; the Poisson bracket of their restrictions to
$\mathbb{S}^2$ can be computed as follows:
\begin{equation}\label{pbs2}
 \{a|_{S^2},b|_{S^2}\}(x)=\langle x,\nabla_xa,\nabla_x b\rangle,
\end{equation}
where $\langle x , \eta , \zeta   \rangle$ stays for  "mixed
product" in $\mathbb{R}^3$, calculated as the determinant of the
$3\times 3$-matrix whose columns are $x , \eta , \zeta$. From now on
we omit the sign of restriction $|_{\mathbb{S}^2}$ while writing
Poisson bracket.

Linear functions $(l ,x )$ are, of course, spherical harmonics. We
denote by $\vec l$ the Hamiltonian field on $\mathbb{S}^2$
associated to the Hamiltonian $\langle l,x\rangle$, $x\in
\mathbb{S}^2$. Obviously, $\vec l$ generates rotation of the sphere
around the axis $l$. According to the aforesaid $\vec l a=\langle
x,l,\nabla a\rangle$ is the Poisson bracket of the functions
$\langle l,x\rangle$ and $a$ restricted to $\mathbb{S}^2$.

The group of rotations acts (by the change of variables) on the
space of homogeneous harmonic polynomials of fixed degree $n$. It is
well-known that this action is irreducible for any $n$ (see
\cite{Arn2} for a sketch of the  proof). In other words, the
following result holds.

\begin{proposition}\label{br1}
Given a nonzero degree $n$ homogeneous harmonic polynomial $a$, the
space
$$
span\{\vec l_1\circ\cdots\circ\vec l_ka:k\ge 0\}
$$
coincides with the space of all degree $n$ homogeneous harmonic
polynomials.
\end{proposition}

\subsection{Poisson brackets of spherical harmonics and
controllability}

Calculating  Lie extensions according to the formula (\ref{fkl0}) we
obtain iterated Poisson brackets of spherical harmonic polynomials,
which in general need not to be harmonic.

The following Lemma shows that there is a way of finding some
harmonic polynomials among them.

\begin{lem}\label{pskob}
  For each $n>2$ there exist  a harmonic homogeneous
polynomial $q$ of degree $2$, and harmonic homogeneous polynomial
$p$ of degree $n>2$ such that their Poisson bracket is again
harmonic (and homogeneous of degree $n+1$) polynomial.
\end{lem}

\begin{proof}
Let us  take so called quadratic zonal harmonic function
$q=\frac{\partial^2 \rho}{\partial x_3^2}$. Being restricted to the
sphere $\mathbb{S}^2$ this function coincides with Legendre
polynomial $q(x_3)=3x_3^2-1$.

Consider homogeneous harmonic polynomials of variables $x_1,x_2$. In
polar coordinates they are known to have representation  $r^{m}\cos
m \varphi$ or alternatively as $Re (x_1+ix_2)^m, \ m=1,2, \ldots $.
We pick the $n$th degree polynomial $p(x_1,x_2)=Re (x_1+ix_2)^n$.

According to (\ref{pbs2}) the Poisson bracket of $q,p$ equals
$$
\{q,p\}= \langle x, \nabla q, \nabla p_n \rangle= \left|
  \begin{array}{ccc}
    x_1 & 0 & p'_{x_1}  \\
    x_2 & 0 & p'_{x_2}  \\
    x_3 & 6x_3 & 0 \\
  \end{array}
\right| =-6x_3\left|
  \begin{array}{ccc}
    x_1 & 0 & p'_{x_1}  \\
    x_2 & 0 & p'_{x_2}  \\
    x_3 & 1 & 0 \\
  \end{array}
\right|.
$$

By  (\ref{pbs2}) the latter determinant  coincides with $\vec e_3
p(x_1,x_2)$, where $e_3=(0,0,1)$ is the standard basis vector of
$\mathbb{R}^3$. Hence by Proposition~\ref{br1} the value of this
determinant is a harmonic polynomial of degree $n$; it equals
$\tilde{p}(x_1,x_2)=-x_1 p'_{x_2}  + x_2p'_{x_1} $ and therefore
does not depend on $x_3$.

Then $\{q,p\}=-6x_3 \tilde{p}(x_1,x_2)$. Since both $-6x_3$ and
$\tilde{p}$ are harmonic, we get $\Delta \{q,p\}=2\nabla(-6x_3)\cdot
\nabla \tilde{p}=-12\partial \tilde{p}/\partial x_3=0$.
\end{proof}

\begin{thm}
Consider NS/Euler system on  sphere $\mathbb{S}^2$. Let (constant)
controlled vector fields correspond to three independent linear
spherical harmonics $l^1,l^2,l^3$, one quadratic harmonic $q$ and
one cubic harmonic $c$. Then this set of controlled vector fields is
saturating and the NS/Euler system is controllable in
finite-dimensional projections and $L_2$-approximately controllable.
$\Box$
\end{thm}

\begin{proof}
It suffices to verify the assumption of the Corollary~\ref{fdens}.
Without lack of generality we may think that $q=\tilde{q}$ - the
second degree zonal harmonic from the previous lemma. Indeed
otherwise we may transform $q$ into $\tilde{q}$  by taking iterated
Poisson brackets with the linear harmonics $l^1,l^2,l^3$.

In fact taking iterated Poisson brackets of $q$ and $c$
respectively, with $l^1,l^2,l^3$ we obtain all quadratic and cubic
harmonics. Thus we manage to obtain all the harmonics of degrees
$\leq 3$.

Let us proceed by induction with respect to the degree of harmonics.
Assume that all harmonics of degrees $\leq n$ are already obtained
by taking iterated Poisson brackets of $\{l^1,l^2,l^3,q,s\}$. Pick
the harmonic polynomial $p$ constructed in Lemma~\ref{pskob}; its
Poisson bracket with $q$ is homogeneous harmonic polynomial
$\bar{p}$ of degree $n+1$. Taking iterated Poisson brackets of
$\bar{p}$ with $l^1,l^2,l^3$ we obtain all polynomials of degree
$n+1$.
\end{proof}

\begin{remark}
Applying the argument similar to the one involved in the pevious
Section one can conclude that there exists a residual set of
Riemannian metrics on $\mathbb{S}^2$ such that the assumptions of
the Corollary~\ref{fdens} are verified and therefore the NS system
is controllable in finite-dimensional projections and
$L_2$-approximately controllable on $\mathbb{S}^2$ endowed with any
of these metrics. $\Box$
\end{remark}

% ----------------------------------------------------------------

\end{document}